\documentclass[11pt]{amsart}

\usepackage{amsfonts,amssymb,stmaryrd,amscd,amsmath,latexsym,amsbsy}

\usepackage{amssymb}
\usepackage{amsfonts}
\usepackage{latexsym}

\newtheorem{theorem}{Theorem}[section]
\newtheorem{lemma}[theorem]{Lemma}
\newtheorem{proposition}[theorem]{Proposition}
\newtheorem{corollary}[theorem]{Corollary}
\theoremstyle{definition}
\newtheorem{definition}[theorem]{Definition}

\newcommand{\g}{\mathfrak{g}}

\newcommand{\ben}{\begin{enumerate}}
\newcommand{\een}{\end{enumerate}}

\newcommand{\mC}{{\mathcal C}}

\newcommand{\CC}{{\mathbb{C}}}

\theoremstyle{plain}

\newtheorem*{sol}{Solution}

\theoremstyle{definition}

\theoremstyle{remark}

\newcommand{\solu}[1]{\begin{sol}{\bf (\ref{#1})}}

\begin{document}

\title{Central extensions of preprojective algebras, the
quantum Heisenberg algebra, and 2-dimensional complex reflection groups}

\author{Pavel Etingof}
\address{Department of Mathematics, Massachusetts Institute of Technology,
Cambridge, MA 02139, USA}
\email{etingof@math.mit.edu}

\author{Eric Rains}
\address{Department of Mathematics, University of California, Davis,
1 Shields Ave, Davis, CA 95616-8633, USA}
\email{rains@math.ucdavis.edu}

\maketitle

\section{Introduction}

Preprojective algebras of quivers were introduced in 1979 
by Gelfand and Ponomarev \cite{GP}, 
because for quivers of finite ADE type, they are models for indecomposable
representations (they contain each indecomposable exactly once). 
Twenty years later, these algebras and their deformed versions introduced 
in \cite{CBH} (for arbitrary quivers) 
became a subject of intense interest, since their
representation varieties, called quiver varieties, 
played an important role in geometric representation theory.  
Ironically, it is exactly for quivers of finite ADE type that 
preprojective algebras fail to have good properties -- they are
not Koszul and their deformed versions are not flat. 

One of the goals of this 
paper is to partially correct this problem. We do so by
introducing a central extension of the preprojective algebra 
of a finite Dynkin quiver (depending on a regular weight
for the corresponding root system), whose natural deformed version 
is actually flat, although it ceases to be flat after
factorization by the central element.\footnote{We note that our
construction makes sense for any quiver, but produces something
new only in the Dynkin case.} We calculate the Hilbert
polynomial of the central extension, and show that it is a
Frobenius algebra. As a corollary, we obtain the Hilbert series
of the usual deformed preprojective algebra in which the
deformation parameters are variables, and show that this algebra
is Gorenstein (although it is not a flat module over the ring of
parameters). 

The main tool in the proofs is 
the fact that our central extension for the weight $\rho$ 
is the image of the quantum Heisenberg algebra in the fusion
category of representations of quantum $SL_2$ under a tensor
functor into $R$-bimodules (where $R$ is the algebra of
idempotents of the quiver). This is a generalization of the result
of \cite{MOV} which says that the usual preprojective algebra is 
the image of the quantum symmetric algebra under the same
functor.  

We also construct Riemann-Hilbert homomorphisms 
from the cyclotomic Hecke algebras of certain 
2-dimensional complex reflection groups to the 
``spherical'' subalgebras of the deformed central extensions
of preprojective algebras. This allows us to show that if all
parameters of the cyclotomic Hecke algebra are equal to 1 (i.e. the
generators are unipotent) then the block of the trivial
representation is equivalent to the category of representations
of the spherical subalgebra of the centrally extended 
preprojective algebra. As a by-product, we show that 
the dimension of the cyclotomic Hecke algebra of a 2-dimensional complex
reflection group for generic parameters is equal to the order of
the group, as conjectured by Broue, Malle, and Rouquier. 

The organization of the paper is as follows. 
In Section 2 we recall basic facts about 
deformed preprojective algebras, Frobenius algebras, and
Cohen-Macaulay and Gorenstein algebras. 
In Section 3 we state the main results regarding central
extensions of preprojective algebras and their deformations,
as well as the corresponding spherical subalgebras. 
In Section 4 we set up the machinery of quantum $SL_2$  
and quantum Heisenberg algebra, which we use in Section 5 to
prove the results of Section 2. Finally, in Section 6 we
introduce and study the Riemann-Hilbert homomorphism.   

{\bf Acknowledgments.} P.E. is grateful to W. Crawley-Boevey
for the references \cite{Ru},\cite{CB2}, and to M. Artin,
J. de Jong, V. Ostrik, D. Rogalski, and R. Rouquier 
for useful discussions. The work of P.E.
was partially supported by the NSF grant DMS-9988796
and the CRDF grant RM1-2545-MO-03.
E.R. was supported in part by NSF Grant No.
DMS-0401387.

\section{Preliminaries} 

\subsection{Quivers and deformed preprojective algebras} 

Let $Q$ be a quiver of finite Dynkin (that is, ADE) type
with Cartan matrix $A$, and Coxeter number $h$.
Let $I=\lbrace{1,...,r\rbrace}$ be the vertex set of $Q$, 
and $R=\oplus_{i\in I}\Bbb C
e_i$ be the commutative algebra generated by idempotents $e_i$
with $e_ie_j=0$, $i\ne j$. Thus any element $\mu\in R$ can be
written as $\mu=\sum_{i\in I}\mu_ie_i$, $\mu_i\in \Bbb C$. 

Let $x_i\in R^*$ be the dual basis to
$e_i$. It is well known that $x_i$ span a root system $\Delta$ inside 
$R^*$ with Cartan matrix $A$.
We will identify $R$ and $R^*$ using the Weyl group invariant inner
product such that $(x_i,x_i)=2$. Then 
$e_i$ are the fundamental weights of $\Delta$.

Let $\rho\in R$ be the sum of fundamental
weights. Of course, $\rho=1$, but we'll use the notation $\rho$
to emphasize the connection with Lie theory. 

Let $\overline Q=Q\cup Q^*$ be the double of $Q$, -- the quiver
with the same vertex set and an additional opposite edge $a^*\in Q^*$ 
for every edge $a\in Q$. 

The generic deformed preprojective algebra $\Pi$ of $Q$
(\cite{CBH},\cite{Ru}) is the 
quotient of the path algebra of $\overline{Q}$
over the ring $\CC[x_1,...,x_r]=S(R^*)$
by the relation 
$$
\sum_{a\in Q}[a,a^*]=\sum_{i\in I} x_ie_i
$$
(here $\sum x_ie_i$ is the canonical element of $R^*\otimes
R$). This algebra is graded ($\text{deg}(e_i)=0$, $\text{deg}(x_i)=2$,
$\text{deg}(a)=\text{deg}(a^*)=1$). Also, it is clear that $\Pi$ is
independent on the orientation of the graph Q, up to an
isomorphism: the reversal of orientation of an edge $a$  
may be accomplished by replacing $a$ by $-a$. 

Let $\Pi_0$ be the zero fiber of $\Pi$, i.e. 
the quotient of the path algebra 
of $\Bbb C\overline{Q}$ of $\overline{Q}$
by the relation $\sum_{a\in Q} [a,a^*]=0$. 
It is called the 
(Gelfand-Ponomarev) preprojective algebra 
of $Q$. The fiber $\Pi_\lambda$ of $\Pi$ at $\lambda\in R$
(i.e. the quotient of the path algebra 
of $\Bbb C\overline{Q}$ of $\overline{Q}$
by the relation $\sum_{a\in Q} [a,a^*]=\lambda$) 
is called the deformed preprojective algebra. 

It is known that $\Pi_0$ is finite dimensional.
Namely, as pointed out already by Gelfand and Ponomarev, 
$\Pi_0$ is a model for indecomposable representations
of $Q$ (i.e., the direct sum of all of them taken once). 
Thus the dimension of $\Pi_0$ is the sum of heights of all
positive roots: 
$$
\dim(\Pi_0)=\sum_{\alpha>0}(\alpha,\rho)=2(\rho,\rho)=\frac{h(h+1)r}{6}.
$$
(see \cite{MOV}). The last equality follows from 
Freudental's magic formula ($(\rho,\rho)=h\dim\g/12$, where $\g$
is the simple Lie algebra attached to $A$) and Kostant's formula 
($(h+1)r=\dim\g$). Moreover, it is known (see \cite{MOV}) that 
the Hilbert polynomial of $\Pi_0$ with respect to its grading is
\begin{equation}\label{h0}
H_0(t)=\frac{1+Pt^h}{1-Ct+t^2}, 
\end{equation}
where $P$ is the permutation of fundamental weights corresponding to
taking the dual representation of $\g$, and $C$ is the adjacency
matrix of $\overline{Q}$.  

This implies that $\Pi_\lambda$ is finite dimensional for each $\lambda$,
and $\Pi$ is a finitely generated $\CC[x_1,...,x_r]$-module.  

It is known \cite{CB1} that unlike the non-Dynkin case, 
the algebra $\Pi$ is not a free $\CC[x_1,...,x_r]$-module. 
More specifically, it is known (see e.g. \cite{CB1}, Theorem 1.2)
that the algebra $\Pi_\lambda$ is zero unless $\lambda$ 
belongs to a reflection hyperplane. 

\subsection{Frobenius algebras}

We recall the basic facts about Frobenius algebras.

Let $A$ be a finite dimensional (unital) algebra over $\CC$. 
Recall that it is called a Frobenius algebra 
if $A\simeq A^*$ as a left $A$-module. 
This is equivalent to saying that there exists 
a linear function $f: A\to \CC$ such that 
the bilinear form $(a,b)=f(ab)$ is nondegenerate. 
Indeed, given $A$ with an isomorphism of left modules 
$\phi: A\to A^*$, we set $f=\phi(1)$, and conversely, 
given $f$, we define $\phi$ by $\phi(a)(b)=f(ba)$. 

\begin{lemma}\label{frob}
Let $A$ be a $\Bbb Z_+$-graded finite dimensional algebra, 
$A=\oplus_{j\ge 0}A[j]$ such that 
$R:=A[0]$ is a commutative semisimple algebra, and the 
Hilbert polynomial $P_A(t)$ satisfies the condition 
$P_A(t)=t^dP_A(t^{-1})$. Then the following conditions are
equivalent:

(i) $A$ is Frobenius;

(ii) $A[d]$ is an invertible $A[0]$-bimodule, and 
the multiplication map \linebreak $A[i]\otimes_R A[d-i]\to A[d]$ 
defines an isomorphism $A[d-i]\to A[i]^*\otimes_R A[d]$. 
\end{lemma}

\begin{proof} Let $R=\oplus_{j=1}^r R_j$, where $R_j$
are copies of $\CC$. Then any $R$-bimodule $M$ can be written as
$M=\oplus M_{ij}\otimes R_{ij}$, where $R_{ij}=\CC$ is the unique
irreducible $(R_i,R_j)$-bimodule. 

Suppose (ii) holds. Since $A[d]$ is invertible, we have 
$A[d]=\oplus_q R_{\sigma(q)q}$ for some permutation $\sigma$. 
Then the multiplication
map defines an isomorphism 
$$
A[d-i]_{pq}\to (A[i]^*)_{p\sigma(q)}.
$$
This implies that if $f\in A^*$ vanishes on degrees $<d$, 
and $f: A[d]\to \CC$ is given by $f(y_1,...,y_r)=\sum_{j=1}^r
y_j$ ($y_p\in R_{p\sigma(p)}$), then $f(ab)$ is a nondegenerate
form on $A$. Thus (i) holds. 

Now suppose (i) holds. Then the map $(a,b)\mapsto f(ab)$ defines 
pairings $A[0]\otimes A[d]\to \CC$ and $A[d]\otimes
A[0]\to \CC$ which have trivial kernel in $A[d]$. 
Since the dimensions of $A[0]$ and $A[d]$ are the same, 
these pairings are nondegenerate. 
This implies that $A[d]$ is faithful as a left and
right $R$-module, which by dimension counting implies that 
$A[d]$ is invertible, i.e. $A[d]=\oplus R_{\sigma(q)q}$. 
Now the multiplication map defines a linear map 
$\psi_{ipq}: A[d-i]_{pq}\to (A[i]^*)_{p\sigma(q)}$.  

Assume that (ii) fails. Then by dimension count 
for some $i,p,q$ the map $\psi_{ipq}$ has a nontrivial kernel. 
Let us choose the largest such $i$ (clearly $i<d$), and let $a$ be an element 
in the kernel of $\psi_{ipq}$. It is clear that $A[m]a=0$ for all
$m>0$, because for any $b\in A[m]_{sp}$, $b\in {\rm Ker}\psi_{i+m,sq}$.
Let us pick $a'\in A[d]$ such that $f(ca)=f(ca')$ for any $c\in
A[0]$. Then $f(c(a-a'))=0$ for all $c\in A$. Contradiction. Thus
(ii) must hold, as desired. 
\end{proof}

\subsection{Cohen-Macaulay and Gorenstein algebras}

Let us now recall the basic properties of noncommutative 
Cohen-Macaulay and Gorenstein algebras. 
(see also \cite{YZ} or the discussion in \cite{EG}, beginning of 
Section 3).

Let $A$ be a $\Bbb Z_+$-graded algebra over $\Bbb C$,
such $R=A[0]$ is finite dimensional and semisimple. 
We will assume that $A$ is a finitely generated module over a finitely 
generated graded central subalgebra $B$, such that $B[0]=\CC$. 
In this case, by Noether's normalization lemma, we may assume
that $B$ is a polynomial algebra; we will do so from now on. 

In this situation, the dualizing complex of $A$ may be
defined by the formula $\Bbb D_A={\rm RHom}_B(A,B)$. 
It can be shown that it does not depend on the choice of $B$. 

\begin{definition} The algebra $A$ is said to be Cohen-Macaulay 
if the cohomology of its dualizing complex 
is concentrated in degree zero, and Gorenstein if this cohomology
is furthermore isomorphic to $A$ as a left module.
\footnote{We note that in the literature there are several
nonequivalent definitions of
Gorenstein property for noncommutative algebras, which are
all equivalent in the commutative case. The definition
we are using implies that $A$ has finite injective dimension, but
is not equivalent to this requirement. We also note that our
definition is not Morita invariant.}
\end{definition} 

Thus we have the following proposition.

\begin{proposition}\label{cm} (i) $A$ is Cohen-Macaulay 
if and only if $A$ is a free $B$-module. 

(ii) $A$ is Gorenstein if and only if furthermore there exists 
a $B$-linear map $f: A\to B$ such that the map $(a,b)\mapsto f(ab)$
is a nondegenerate pairing $A\times A\to B$, i.e. defines an
isomorphism $A\to {\rm Hom}_B(A,B)$. 
\end{proposition}    

This immediately implies the following. 

\begin{proposition}\label{gore}  
(i) Let $A$ be a finite dimensional algebra. 
Then $A$ is Gorenstein if and only if it is Frobenius. 

(ii) $A$ is Gorenstein if and only if its zero-fiber $A_0$ (as a
$B$-module) is Gorenstein.\footnote{Here by definition 
$A_0=A/JA$, where $J$ is the kernel of the augmentation in $B$.} 
\end{proposition}

\section{Central extensions of preprojective algebras}

\subsection{Gorenstein properties and Hilbert series 
of generic deformed preprojective algebras}

One of the main results of this paper is the following theorem.

\begin{theorem}\label{gor1}
(i) $\Pi$ is a Gorenstein algebra (in particular, Cohen-Macaulay).

(ii) The matrix Hilbert series of $\Pi$ 
(i.e. the matrix consisting of the Hilbert series 
of $e_i\Pi e_j$) is equal to 
$$
H(t)=\frac{1-t^{2h}}{(1-t^2)^r(1-Ct+t^2)},
$$
where $C=2-A$ is the adjacency matrix of $\overline Q$. 
\end{theorem}

Theorem \ref{gor1} is proved in Section \ref{Proofs}.

\subsection{Central extensions of preprojective algebras}

Another result of this paper is a construction of a central
extension of the preprojective algebra, whose deformed version 
is flat, unlike that of the usual preprojective algebra. 
Namely, let $\Pi[z]$ be the algebra 
of polynomials of a central variable $z$ with coefficients in
$\Pi$; it is graded with $\text{deg}(z)=2$. 
Let $\mu\in R$ be a regular weight (i.e., does not belong to any
reflection hyperplane $(\alpha,\mu)=0$).
For any weight $\lambda\in R$, let $\Pi_\lambda^\mu$ be
the quotient of $\Pi[z]$ by the relations 
$x_i=\mu_iz+\lambda_i$. Thus $\Pi_\lambda^\mu$ is the quotient of the 
path algebra $\CC[z]\overline{Q}$
by the relation 
$$
\sum_a [a,a^*]=\sum_i (\mu_iz+\lambda_i)e_i. 
$$
This algebra carries a natural filtration induced by the grading
in the path algebra. 

\begin{theorem}\label{cenex}
(i) $\Pi[z]$ is a finitely generated free module over 
the algebra $K_\mu:=\CC[x_1-\mu_1z,...,x_r-\mu_rz]$.

(ii)  The associated graded
algebra of $\Pi_\lambda^\mu$ under the natural filtration is 
$\Pi_0^\mu$. 

(iii) The algebra $\Pi_0^\mu$ is Gorenstein (=Frobenius). 
One has $z^{h-1}=0$, and for generic $\mu$ 
the socle of $\Pi_0^\mu$ is $z^{h-2}R$.

(iv) The Hilbert polynomial of $\Pi_0^\mu$ 
is 
$$
\widetilde
H(t)=H(t)(1-t^2)^{r-1}=\frac{1+t^2+...+t^{2(h-1)}}{1-Ct+t^2}.
$$

(v) The dimension of $\Pi_\lambda^\mu$ is $h^2(h+1)r/12$,
and for generic $\lambda$ it is semisimple. 
\end{theorem}

Theorem \ref{cenex} is proved in Section \ref{Proofs}.

Part (v) of Theorem \ref{cenex} implies that for generic
$\lambda$ the element $z$ is semisimple, and hence 
$$
\Pi_\lambda^\mu=\oplus_{\alpha>0}\Pi_{\lambda-\frac{(\lambda,\alpha)}
{(\mu,\alpha)}\mu}.
$$
Thus $\Pi_\lambda$ is semisimple for a generic $\lambda$ on a
reflection hyperplane $(\alpha,\lambda)=0$ (where $\alpha$ is a
positive root). 
Moreover, by \cite{CB1}, Theorem 1.2, it is actually
simple: $\Pi_\lambda={\rm Mat}_{(\alpha,\rho)}(\Bbb C)$.
This means that the rational Weyl denominator
$\delta(x):=\prod_{\alpha>0}
(\alpha,x)\in \CC[x_1,...,x_r]$ 
is zero in $\Pi$, so $\Pi$ is scheme-theoretically supported 
on the reflection hyperplanes. This was conjectured by Rump
\cite{Ru} and proved by Crawley-Boevey in \cite{CB2} using a
different method.  

\subsection{Subalgebras corresponding to nodal vertices}\label{nv}

It is interesting to consider ``spherical subalgebras'' of the above
algebras corresponding to nodal vertices. 

Namely, let $Q$ be of type $A_{2n-1}$, $D$ or $E$, and 
$p\in I$ be the nodal vertex (= the branching vertex for $D$ and $E$, and
the middle vertex for $A_{2n-1}$).  
After removal of $p$, the quiver becomes a union of ``legs'',
i.e. Dynkin diagrams of type $A_{d_k-1}$, $k=1,...,m$ ($m=2$ or $3$). 
Recall that to such $Q$ one can attach a finite 
subgroup $G$ of $SO(3)$ generated by elements $T_k$, $k=1,...,m$, 
with defining relations $T_k^{d_k}=1$ for $k=1,...,m$, and 
$\prod_{k=1}^m T_k=1$ (for type $A$ we get the cyclic group, for type $D$ the dihedral group, 
for type $E_6$ the tetrahedral group, for type $E_7$ the cube
group, and for type $E_8$ the icosahedral group); all finite
subgroups of $SO(3)$ are obtained in this way. 

The order $|G|$ of the group $G$ can be represented as a product 
of two integers $q_1\le q_2$ such that $q_1+q_2-1=h/2$. 
Namely (see e.g. \cite{Ko}), $q_1$ is the p-th coordinate of the maximal root 
in the basis of simple roots.
\footnote{Another context in which the numbers $q_1,q_2$ appear
is the following: the Hilbert series of the invariants
$\CC[x,y]^\Gamma$ 
in $\Bbb C[x,y]$ under the double cover $\Gamma\subset
SL(2)$ of $G$ is $\frac{1+t^h}{(1-t^{2q_1})(1-t^{2q_2})}$.}

Let $B=e_p\Pi e_p$. 
According to \cite{MOV}, \cite{Me} the algebra $B$ is generated 
over $\Bbb C[x_1,...,x_r]$ by $U_k$, 
$k=1,...,m$, with defining relations 
$$
U_k(U_k-x_{i_1(k)})...(U_k-x_{i_1(k)}-...-x_{i_{d_k-1}(k)})=0,
$$
where $i_1(k),...,i_{d_k-1}(k)$ are the vertices
of the $k$-th leg of $Q$ enumerated from the nodal vertex, and
$$
\sum_{k=1}^m U_k=-x_p.
$$
Namely, the elements $U_k$ are just the elements
$a_k^*a_k$, where $a_k$ are the edges of $\overline Q$ starting
at $p$ and going
along the $k$-th leg.

This implies that the algebra
$B^\mu_\lambda:=e_p\Pi^\mu_\lambda e_p$
is generated over $\Bbb C[z]$ by $U_k$
with the same defining relations, in which 
$x_i=\lambda_i+\mu_iz$. \footnote{Note that in our normalization,
the degree of all generators of this algebra, including $z$, is $2$.} 

Let $\mu$ be a weight such that $(\alpha,\mu)\ne 0$
for any positive root $\alpha$ that involves $x_p$ 
with a strictly positive coefficient. 
Also, for a number $n$, let $[n]_q:=\frac{1-q^n}{1-q}$.
For brevity, if $q=t^2$, we will simply write $[n]$ instead of $[n]_q$.

\begin{theorem}\label{Bz} 
(i) The algebra $B[z]$ is a finitely generated free module over 
the algebra $K_\mu$.

(ii)  The associated graded
algebra of $B_\lambda^\mu$ under the natural filtration is 
$B_0^\mu$. 

(iii) The algebra $B_0^\mu$ is Gorenstein. We have $z^{h-1}=0$,
and for generic $\mu$ the socle of $B_0^\mu$ is spanned by $z^{h-2}$. 

(iv) The Hilbert polynomial of $B_0^\mu$ 
is $E(t)=[h/2][q_1][q_2]$.

(v) The dimension of $B_\lambda^\mu$ is $hq_1q_2/2$, and for
generic $\lambda$ it is semisimple. 
\end{theorem}

Theorem \ref{Bz} is proved in Section \ref{Proofs}.

{\bf Remark.} If $\mu$ is regular, then this theorem easily
follows from Theorem \ref{cenex}, but our statement is 
more general, so it requires a separate proof. 

Considering the case $\mu=e_p$, 
we immediately get 
the following corollary. For any 
numbers $\lambda_{ik}$, $i=1,...,d_k$, $k=1,...,m$, define the
algebra $B(\lambda)$ generated by $U_k$ and a central element $z$ 
with defining relations 
$$
\prod_{i=1}^{d_k} (U_k-\lambda_i)=0,
$$
$$
\sum_{k=1}^m U_k=z.
$$

\begin{corollary} \label{defo}
(i) ${\rm gr}B(\lambda)=B(0)$.

(ii) $B(0)$ is a Gorenstein algebra, and $z^{h-1}=0$ in $B(0)$.

(iii) The Hilbert polynomial of $B(0)$ is $[h/2][q_1][q_2]$.
\end{corollary}

Indeed, $B(\lambda)$ is obtained from $B_{\lambda'}^{e_p}$ 
by a straightforward change of variables. 

For comparison note that the Hilbert polynomial of the algebra 
$B_0:=B(0)/(z)$, according to \cite{MOV}, equals $[q_1][q_2]$. 

{\bf Remark.} We have checked using the Magma computer algebra
system that the socle of $B(0)$ is spanned by $z^{h-2}$.

\section{The quantum Heisenberg algebra}

\subsection{Definition and properties of the quantum Heisenberg
algebra}

The proofs of the main results of this paper are based on 
the idea of \cite{MOV}: algebras related to quivers 
may be obtained from algebras in the category of representations of
quantum $SL(2)$ by application of tensor functors into
$R$-bimodules.

We are going to use the notions and notation from \cite{BaKi,EO}. 
Namely, for every $q\in \Bbb C^*$, denote by 
$\widetilde \mC_q$ the tensor category of finite dimensional comodules
over the quantum function algebra $F_q(SL(2))$. If $q$ is not a
root of unity of order $>2$, then this category is semisimple and
has simple objects $V_i$, $i\in \Bbb Z_+$ (representations with highest
weight $i$) such that $V_0$ is the neutral object, 
with the Clebsch-Gordan tensor product rule 
$$
V_i\otimes V_j=\bigoplus_{n=0}^{\min(i,j)}V_{2n+|i-j|}. 
$$
On the other hand, if $q$ is a root of unity of order $n>2$, then the category 
$\widetilde \mC_q$ is not semisimple. However, the objects 
$V_i$ are still well defined (the so called dual Weyl modules,
i.e. the homogeneous components of the quantum symmetric algebra), and 
they are simple if $i\le \bar n-1$, where $\bar n=n$ if $n$ is
odd and $\bar n=n/2$ if $n$ is even. 

Consider the tensor
algebra $T(V_1\oplus V_0)$ in $\widetilde 
\mC_q$, an ind-object in $\widetilde 
\mC_q$. 
We will regard it as a graded algebra in which 
$V_1$ has degree $1$ and $V_0$ has degree 2. 
Then we have two morphisms from $V_0$ 
to $T(V_1\oplus V_0)[2]$: the map 
$f$ which is a composition 
$V_0\to V_1\otimes V_1\to T(V_1\oplus V_0)[2]$
\footnote{Here the first map is a generator of
the space ${\rm Hom}(V_0,V_1\otimes V_1)$, and the second map is
induced by multiplication.}, and 
the map $g: V_0\to T(V_1\otimes V_0)[2]$ coming from the
embedding $V_0\to V_1\oplus V_0$. 
There are also two obvious embeddings 
$f_1,f_2: V_0\otimes V_1 \to T(V_0\oplus V_1)[3]$ 
corresponding to multiplication in two different orders 
(more precisely, 
$f_1$ is the multiplication map, and 
$f_2$ is the multiplication map composed with the canonical isomorphism 
$V_0\otimes V_1\to V_1\to V_1\otimes V_0$). 
Let $\widetilde{\mathcal J}$ be the ideal 
generated by the images of $f-g$ and $f_1-f_2$. 

\begin{definition}
The quantum Heisenberg algebra in $\widetilde \mC_q$ 
is the algebra $\widetilde A=T(V_0\oplus
V_1)/\widetilde{\mathcal J}$. 
\end{definition}

More explicitly, 
the algebra $\widetilde A$ is the usual quantum Heisenberg
algebra generated by $x,y,z$ with defining relations 
saying that the element $z=xy-qyx$ is central\footnote{Here $x,y$
is a weight basis of $V_1$, and $z$ is a basis element of $V_0$.}, i.e.
$$
xy-qyx=z,\ xz-zx=0,\ yz-zy=0,
$$
and the usual coaction of the Hopf algebra $F_q(SL(2))$;
it is convenient to express this coaction as an action 
of $U_q(sl(2))$ given by 
$$
ex=y,\ ey=0,\ fy=x,\ fx=0,\ q^hy=qy,\ q^hx=q^{-1}x,\ 
$$
$$
ez=fz=0,\ q^hz=z.  
$$
(If $q$ is a root of unity, one should add that the divided
powers of $e$ and $f$ act by zero). 
Here the coproduct of $U_q(sl(2))$ is given by 
$\Delta(e)=e\otimes q^h+1\otimes e$, $\Delta(f)=f\otimes
1+q^{-h}\otimes f$, $\Delta(q^h)=q^h\otimes q^h$. 

It is not hard to show that the elements $y^ix^jz^m$, $i,j,m\ge 0$, 
form a basis of $\widetilde A$, and the multiplication in this
basis is derived from the commutation relation
$$
x^py^j=\sum_{i=0}^pq^{(j-i)(p-i)}\binom{p}{i}_q\prod_{s=1}^i[j-s+1]_q\cdot
y^{j-i}x^{p-i}z^i,
$$
where $\binom{p}{i}_q:=[p]_q!/[i]_q![p-i]_q!$, and
$[p]_q!:=[1]_q...[p]_q$.

Now let us consider the structure of the homogeneous subspaces
$\widetilde{A}[j]$ as $F_q(SL(2))$-comodules. 

\begin{lemma}\label{nondeg}
(i) If $q$ is not a root of unity (of order $>1$) 
then, as an $F_q(SL(2))$-comodule, $\widetilde A[j]$ 
is isomorphic to the direct sum $V_j\oplus V_{j-2}\oplus...\oplus
V_{\bar j}$, where $\bar j$ 
is $j$ modulo $2$. 

(ii) If $q$ is a root of unity and the order of $q$ is $2h$, where
$h\ge 2$, then the decomposition of (i) holds for $j<h$. 
Furthermore, for $j=h+p-1$ with $h-1\ge p\ge 1$, the representation 
$z^p\widetilde{A}[j-2p]=V_{h-p-1}\oplus V_{h-p-3}\oplus...$
is canonically contained in $\widetilde A[j]$. 

(iii) In the situation of (ii) let $j<h$, and $E_{i,s}$ be the
canonical copy of $V_{i-2s}$ in $\widetilde A[i]$. Then the image $Y_j$
of the multiplication map 
$E_{j,0}\otimes E_{1,0}  \to \widetilde{A}[j]$ contains $E_{j+1,1}$.
\end{lemma}

\begin{proof}
Parts (i) and (ii) are straightforward by looking at the 
characters of the representations. 

Let us prove part (iii). It is easy to see that

\begin{equation}\label{commu}
xy^j=q^jy^jx+[j]_qy^{j-1}z.
\end{equation}
On the other hand, 
$$
fe(y^jx)=f(y^{j+1})=\sum_{s=0}^j q^{-s} y^sxy^{j-s}.
$$
Therefore, using equation (\ref{commu}), 
we get 
$$
fe(y^jx)=\sum_{s=0}^j q^{-s}y^s(q^{j-s}y^{j-s}x+[j-s]_qy^{j-s-1}z),
$$
which yields
\begin{equation}\label{fe}
fe(y^jx)=\frac{q^{j+1}-q^{-j-1}}
{q-q^{-1}}y^jx+\frac{[j]_q[j+1]_{q^{-1}}}{[2]_{q^{-1}}}y^{j-1}z.
\end{equation}
Since the second coefficient is not zero, 
and $y^jx\in Y_j$, we get that $y^{j-1}z\in Y_j$
(as $Y_j$ is a subrepresentation of the quantum group). 
But $y^{j-1}z$ is the highest weight vector of $E_{j+1,1}$,
which yields (iii). 
\end{proof} 

Now let $q=e^{\pi i/h}$, where $h\ge 2$ is an integer. Then  
the category $\widetilde \mC_q$ 
contains a non-abelian monoidal subcategory ${\mathcal T}_q$
of tilting modules, which is closed under taking direct summands.
Furthermore, the category ${\mathcal T}_q$ contains a tensor ideal
${\mathcal I}$, such that ${\mathcal T}_q/{\mathcal I}=\mC_q$ 
is a semisimple tensor category (the fusion category). 
The simple objects in $\mC_q$ are $V_0,...,V_{h-2}$, 
and the tensor product is given by the Verlinde 
rule: 
$$
V_i\otimes V_j=\bigoplus_{n=0}^{\min(i,j,h-2-\max(i,j))}V_{2n+|i-j|}. 
$$

Consider now the tensor algebra $T(V_0\oplus V_1)$ in the fusion
category $\mC_q$. Then we can define the ideal $\mathcal J$ 
in this algebra generated by the images of the maps $f-g$ and
$f_1-f_2$, where $f_1,f_2,f,g$ are defined in the same way as 
above. 

\begin{definition}
The quantum Heisenberg algebra in $\mC_q$ 
is the algebra $A=T(V_0\oplus
V_1)/\mathcal J$. 
\end{definition}

{\bf Remark.} The quotient $A_0$ of $A$ by the ideal generated by
the unique copy of
$V_0$ in degree $2$ is the quantum symmetric algebra considered
in \cite{MOV}; as an object of $\mC_q$ it is
$V_0\oplus V_1\oplus...\oplus V_{h-2}$, and the degree 
of $V_i$ is $i$. 

\begin{proposition}\label{alg}
(i) The degree $n$ component $A[n]$ of $A$ is the object 
$\oplus_{j\le s/2}V_{s-2j}$, where $s={\rm min}(n,2h-4-n)$. 

(ii) The algebra $A$ is Gorenstein in the following sense: 
the highest nontrivial degree in $A$ is $A[2h-4]=V_0$, and 
the multiplication map \linebreak $A[i]\otimes A[2h-4-i]\to A[2h-4]$ defines
an isomorphism $A[i]\to A[2h-4-i]^*$. 
\end{proposition}

\begin{proof}
Let $L$ be the ideal in $\widetilde A$ generated by the copy 
of $V_{h-1}\subset V_1^{\otimes h-1}\subset T(V_0\oplus
V_1)[h-1]$. It is obvious that for any $j=0,...,h-1$, 
$L[h-1+j]$ has trivial intersection with the canonical
subrepresentation $V_{h-3-j}\oplus V_{h-5-j}\oplus...$.  
On the other hand, it follows from Lemma \ref{nondeg}(iii) that 
$L[h-1+j]$ contains the canonical copy of $V_{h-1-j}$
for all $0\le j\le h-1$. 

By looking at the associated graded
algebra of $\widetilde A$ under the filtration defined 
by $\deg x=\deg y=1$, $\deg z=0$, this implies that 
$L[h-1+j]$ is a complement of $V_{h-3-j}\oplus
V_{h-5-j}\oplus...$ in $\widetilde{A}[h-1+j]$. 

Therefore, the graded algebra 
$\widetilde A/L$, as an object of the
category $\widetilde \mC_q$, has the structure specified in 
part (i) of the proposition. In particular, $\widetilde A/L$ belongs 
to the category of tilting modules $\mathcal T_q$.
In the category $\mathcal T_q$, the algebra $\widetilde A/L$ can
be written as $T(V_0\oplus V_1)/\widetilde L$, where 
$\widetilde L$ is the preimage of $L$ in $T(V_0\oplus V_1)$. 
Note that both the algebra $T(V_0\oplus V_1)$ and the ideal
$\widetilde L$ belong to $\mathcal T_q$. 

Now we can consider the image $A'$ of $\widetilde A/L$ in the fusion category
$\mC_q$ (which is by definition a quotient category of $\mathcal
T_q$). We have $A'=T(V_0\oplus V_1)/J'$, where $J'$ is the image
of $\widetilde L$. Since the image of $V_{h-1}$ in $\mC_q$ is
zero, we have $J'=J$ and $A'=A$, hence (i). 
 
To prove part (ii), it suffices to note that by Lemma
\ref{nondeg}(iii), if $V_j\subset A[i]$, $i\ge j>0$ then 
the map $V_j\otimes V_1\to V_{j-1}$ defined by 
multiplication by the generating copy of $V_1$ and then
projecting to $V_{j-1}\subset A[i+1]$ is nonzero.
More precisely, this statement follows from the lemma by
multiplying by $z^{(i-j)/2}$.  
\end{proof}

\subsection{The quantum Heisenberg algebra and the central
extension of the preprojective algebra}

Now assume that $q=e^{\pi i/h}$, where $h$ is the Coxeter number of
$Q$. Recall (see \cite{EO,MOV}) that there exists a unique tensor
functor $\mathcal F: \mC_q\to R-bimod$ from $\mC_q$ to the category of
$R$-bimodules, such that the bimodule $\mathcal F(V_1)$ is the edge space
$E$ of the doubled quiver $\overline{Q}$. 
It is checked in \cite{MOV} that $\mathcal F(A_0)$ is the 
preprojective algebra $\Pi_0$. 

\begin{proposition}\label{func}
$\mathcal F(A)=\Pi_0^\rho$.
\end{proposition}

\begin{proof}
The algebra $\Pi_0^\rho$ is the quotient of the path algebra 
$\Bbb C[z]\overline Q$ 
by the relations saying that $\sum_{a\in Q}[a,a^*]=z$, and $z$ is
central. It is easy to show, similarly to \cite{MOV}, that 
these relations are images under $\mathcal F$ of the relations $xy-qyx=z$, $xz=zx$,
$yz=zy$. The proposition is proved. 
\end{proof}

\begin{corollary}\label{pir} The Hilbert series of $\Pi_0^\rho$
is given by the formula in Theorem \ref{cenex}, (iv). In
particular, the dimension of $\Pi_0^\rho$ is 
$h^2(h+1)r/12$.  
\end{corollary} 

\begin{proof} In the Grothendieck ring of 
$\mC_q$, we have $V_j=P_j(V_1)$, where 
$P_j$ is the Tchebysheff polynomial 
of the second kind: $P_j(2\cos x)=\frac{\sin (j+1)x}{\sin x}$.
This means that by Proposition \ref{alg}, (i), 
the Grothendieck-group-valued Hilbert polynomial of $A$ is 
$$
\sum_{j=0}^{h-2}\sum_{i=0}^{h-2-j}t^{2i+j}P_j(V_1), 
$$
and hence by Proposition \ref{func}, 
the Hilbert polynomial of $\Pi_0^\rho$ is 
$$
\widetilde H(t)=\sum_{j=0}^{h-2}\sum_{i=0}^{h-2-j}t^{2i+j}P_j(C)=
\sum_{j=0}^{h-2}\frac{t^j-t^{2(h-1)-j}}{1-t^2}P_j(C)=
\frac{H_0(t)-t^{2h-2}H_0(t^{-1})}{1-t^2}.
$$
Thus by formula (\ref{h0}),
$$
\widetilde
H(t)=\frac{1-t^{2h}}{(1-t^2)(1-Ct+t^2)},
$$
as desired. Thus we have $\widetilde
H(t)=\frac{1-Pt^h}{1-t^2}H_0(t)$, and hence 
$\dim \Pi_0^\rho=h\dim\Pi_0/2=h^2(h+1)r/12$. 
The corollary is proved. 
\end{proof}

{\bf Remark.} Since by \cite{MOV}, ${\mathcal F}(A_0)=\Pi_0$, we 
conclude from Proposition \ref{alg}
that the algebra $\Pi_0$, and hence $B_0$, are
Frobenius algebras. 

\section{Proofs of Theorems \ref{gor1},\ref{cenex},\ref{Bz}}\label{Proofs}

\subsection{Proofs of Theorems \ref{gor1},\ref{cenex}}

We now prove Theorems \ref{gor1},\ref{cenex}. 

\begin{lemma}\label{pir1} $\Pi[z]$ is a free module 
over $\CC[x_1-z,...,x_r-z]$. 
\end{lemma}

\begin{proof} The fiber of this module at $x_i-z=\lambda_i$ is
$\Pi_\lambda^\rho$. Thus, since $\Pi[z]$ is $\Bbb Z_+$-graded, 
it suffices to show
that the dimension of $\Pi_\lambda^\rho$ is at least as big as
that of $\Pi_0^\rho$, for generic $\lambda$. To do so, 
note that by Theorem 1.2 in \cite{CB1}, for generic $\lambda$ 
the algebra $\Pi_\lambda^\rho$ has an irreducible
representation with dimension vector being each positive root
$\alpha$. Therefore, 
$$
\dim \Pi_\lambda^\rho\ge \sum_{\alpha>0}(\alpha,\rho)^2=
h(\rho,\rho)=h^2\frac{\dim \g}{12}=\frac{h^2(h+1)r}{12},
$$
which together with Corollary \ref{pir} 
implies the desired inequality.
\end{proof} 

Lemma \ref{pir1} together with Corollary \ref{pir} 
imply part (ii) of Theorem \ref{gor1}. 

Now, by Proposition \ref{alg}, Proposition \ref{func}, and
Lemma \ref{frob}, 
$\Pi_0^\rho$ is a Frobenius algebra. Hence by Proposition \ref{gore} (i) 
it is a Gorenstein algebra. 
Therefore, since $\Pi[z]$ is a free module 
over $\CC[x_1-z,...,x_r-z]$, by Proposition \ref{gore} (ii),
the algebra $\Pi[z]$ is Gorenstein, and hence 
again by Proposition \ref{gore} (ii), $\Pi$ is Gorenstein. 
This implies part (i) of Theorem \ref{gor1}. 
Thus Theorem \ref{gor1} is proved.  

Now let $\mu$ be a regular weight. Since $\Pi_0$ is finite
dimensional, the algebra $\Pi_0^\mu$ is finitely generated as
a module over $\Bbb C[z]$. Moreover, if $z_0\in \CC$, we have 
$\Pi_0^\mu/(z-z_0)=\Pi_{z_0\mu}$, which is zero by Theorem 1.2 of
\cite{CB1} for any $z_0\ne 0$. Hence $\Pi_0^\mu$ is finite
dimensional. This implies that $\Pi[z]$ is a finitely generated
$K_\mu$-module (since this module is $\Bbb Z_+$-graded and its 
zero-fiber is $\Pi_0^\mu$). Since $\Pi[z]$ is a Gorenstein
algebra, it is in particular Cohen-Macaulay, and by 
Proposition \ref{cm} (i), $\Pi[z]$ is a free $K_\mu$-module,
hence parts (i),(ii) of Theorem \ref{cenex}. 

We also see that the Hilbert polynomial of $\Pi_0^\mu$ is the
same as that for $\Pi_0^\rho$, so 
Corollary \ref{pir} implies Theorem \ref{cenex}, (iv), and the
dimension formula of (v). The semisimplicity of 
$\Pi_\lambda^\mu$ follows from the fact that for generic
$\lambda$, by Theorem 1.2 of \cite{CB1}, it has an irreducible
representation with dimension vector being every positive root
$\alpha$; as we've shown, the sum of dimensions of these
representations is $h^2(h+1)r/12=\dim \Pi_0^\mu$. 

We also see, by Proposition \ref{gore} (ii), that $\Pi_0^\mu$ is
Gorenstein (=Frobenius). This proves the first statement of 
Theorem \ref{cenex}, (iii). The second statement follows
from part (iv), and the third statement (the socle is $z^{h-2}R$)
follows from the fact that this is so for $\mu=\rho$ (by
Proposition \ref{alg} and Proposition \ref{func}). 
Theorem \ref{cenex} is proved. 

{\bf Remark.} Let $N$ be the ideal in $\Pi_0^\mu$ 
generated by $z$, and $H_k(t)$ be the Hilbert polynomial 
of the quotient $N^k/N^{k+1}$. The above arguments show that for
generic $\mu$ one has 
$$
\sum_{k\ge 0}H_k(t)u^k=\frac{H_0(t)-ut^hPH_0(ut)}{1-ut^2}.
$$
Indeed, if $\mu=\rho$ then using Propositions \ref{alg} and \ref{func} as in the
proof of Corollary \ref{pir}, we get 
$$
\sum_{k\ge 0}H_k(t)u^k=
\sum_{j=0}^{h-2}\sum_{i=0}^{h-2-j}t^{2i+j}u^jP_j(C)=
$$
$$
\sum_{j=0}^{h-2}\frac{t^j-t^{2(h-1)-j}u^{h-1-j}}{1-ut^2}P_j(C)=
\frac{H_0(t)-t^{2h-2}u^{h-1}H_0(t^{-1}u^{-1})}{1-ut^2},
$$ 
which gives the desired formula. Also, it is clear from this
argument that for $\mu=\rho$ the
operator $z: \Pi_0^\mu[j]\to \Pi_0^\mu[j+2]$ has maximal rank for
all $j$. Hence, the same formula applies to generic $\mu$. 

\subsection{Proof of Theorem \ref{Bz}}

It follows from Theorems \ref{gor1},\ref{cenex} that
the algebras $B$ and $B[z]$ are Gorenstein. 
Since $B_0$ is finite dimensional, the algebra $B_0^\mu$ is finitely generated as
a module over $\Bbb C[z]$. Moreover, by Theorem 1.2 of
\cite{CB1}, $z$ must act by zero in every irreducible
representation of $B_0^\mu$. Thus, $B_0^\mu$ is finite
dimensional. This implies that $B[z]$ is a finitely generated
$K_\mu$-module. Since $B[z]$ is a Gorenstein
algebra, it is in particular Cohen-Macaulay, and by 
Proposition \ref{cm} (i), $B[z]$ is a free $K_\mu$-module,
hence parts (i),(ii) of Theorem \ref{Bz}. 

We also see that the Hilbert polynomial 
of $B_0^\mu$ is the
same as that for $B_0^\rho$, so 
Proposition \ref{pir} implies Theorem \ref{Bz}, (iv), and the
dimension formula of (v). The semisimplicity of 
$B_\lambda^\mu$ follows from the fact that for generic
$\lambda$, by Theorem 1.2 of \cite{CB1}, it has an irreducible
representation corresponding to every positive root
$\alpha$ that involves $x_p$; the sum of dimensions of these
representations is $hq_1q_2/2=\dim B_0^\mu$. 

We also see, by Proposition \ref{gore} (ii), that $B_0^\mu$ is
Gorenstein (=Frobenius). This proves the first statement of 
Theorem \ref{Bz}, (iii). The second statement follows
from part (iv), and the third statement (the socle is $z^{h-2}$)
follows from the fact that this is so for $\mu=\rho$ (by
Proposition \ref{alg} and Proposition \ref{func}). 
Theorem \ref{Bz} is proved. 

\section{Relation to Hecke algebras of 
2-dimensional complex reflection groups}

Let us discuss the connection of the above results with the theory of
Hecke algebras for complex reflection groups, due to 
Brou\'e, Malle, and Rouquier (\cite{BMR}).

\subsection{Finiteness of cyclotomic Hecke algebras in 2 dimensions}

Let $\Gamma$ be a 2-dimensional irreducible complex reflection group, and 
$H$ be its cyclotomic Hecke algebra (see \cite{BMR}); it is
defined by generators and relations which are deformations 
of the relations of $\Gamma$. The algebra
$H$ is a module over the algebra $\Bbb C[\bold T]$ 
of functions on the torus $\bold T$ of parameters. Denote by
$H(\Lambda)$ 
the specialization of $H$ at a point $\Lambda\in \bold T$. 

It is known (\cite{BMR}) that the dimension
of $H(\Lambda)$ is generically $\ge |\Gamma|$.
Moreover, it is conjectured (see \cite{BMR}), that the algebra 
$H$ is a free $\Bbb C[\bold T]$-module, 
and hence $H(\Lambda)$ has dimension exactly 
$|\Gamma|$ for all values of $\Lambda$, i.e.  
it is a flat deformation of the group
algebra $\Bbb C[\Gamma]$. As far as we know, 
this conjecture has been
checked by J. M\"uller using a computer 
for all cases except $G_{17},G_{18},G_{19}$
(see \cite{GGOR}, Remark 5.12).

Here we give a proof 
of a weak version of this conjecture. 

\begin{theorem}\label{finitel}
The algebra $H$ is a finitely generated module over $\CC[\bold
T]$. Therefore, the algebra $H(\Lambda)$ is finite
dimensional for all $\Lambda$, and has dimension $|\Gamma|$ for generic
$\Lambda$, and $\ge |\Gamma|$ for all $\Lambda$. 
\end{theorem}

The rest of the subsection is devoted to the proof of this
theorem\footnote{The first author is very grateful to R. Rouquier
for help with this proof.}. 

The theorem is known for the infinite series 
of complex reflection groups (\cite{BMR}), 
so it is sufficient to concentrate on 
the exceptional ones, i.e. the groups $G_n$ in the Shephard-Todd
classification \cite{ST}, $4\le n\le 22$.

For the proof we need explicit presentations of the 
groups $\Gamma$. Recall that in all cases the image $G$ of $\Gamma$ in
$PGL_2(\CC)$ is the tetrahedral group, 
octahedral group, or icosahedral group. 
Thus the groups $G_4,...,G_{22}$ fall into 
three families: tetrahedral, octahedral, and icosahedral. 
In each family, there is a maximal group $\widehat G$ 
of order $|G|^2$, which is a central extension of $G$ by a cyclic
group of order $|G|$ (generated by an element $Z$), 
and all others are its subgroups. 

The realizations of the particular groups $G_n$, $4\le n\le 22$, 
are as follows. 

{\bf The tetrahedral family.}

The tetrahedral group $G$: $a^2=b^3=c^3=1$, $abc=1$ (order 12).

The maximal group $\widehat G=G_7$ (order 144):
$a_*^2=b_*^3=c_*^3=1$, $a_*b_*c_*=Z$, $Z$ is central; conjugacy classes
of reflections are powers of $a_*,b_*,c_*$. 

Group $G_4$ (order 24) is generated by $a=a_*Z^{-3},b=b_*Z^2,c=c_*$. 
The relations are $a^2=\zeta^{-1}$, $b^3=\zeta$, $c^3=1$,
$abc=1$, $\zeta$ is central (here $\zeta=Z^6$). Conjugacy classes
of reflections are powers of $c$. 

Group $G_5$ (order 72) is generated by $a=a_*Z^{-1},b=b_*,c=c_*$. 
The relations are $a^2=\zeta^{-1}$, $b^3=1$, $c^3=1$,
$abc=1$, $\zeta$ is central (here $\zeta=Z^2$). Conjugacy classes
of reflections are powers of $b,c$. 

Group $G_6$ (order 48) is generated by $a=a_*,b=b_*Z^{-1},c=c_*$. 
The relations are $a^2=1$, $b^3=\zeta^{-1}$, $c^3=1$,
$abc=1$, $\zeta$ is central (here $\zeta=Z^3$). Conjugacy classes
of reflections are powers of $a,c$. 

{\bf Octahedral family.}

The octahedral group $G$: $a^2=b^3=c^4=1$, $abc=1$ (order 24).

The maximal group $\widehat G=G_{11}$ (order 576):
$a_*^2=b_*^3=c_*^4=1$, $a_*b_*c_*=Z$, $Z$ is central; conjugacy classes
of reflections are powers of $a_*,b_*,c_*$. 

Group $G_8$ (order 96) is generated by $a=a_*Z^{-3},b=b_*Z^2,c=c_*$. 
The relations are $a^2=\zeta^{-1}$, $b^3=\zeta$, $c^4=1$,
$abc=1$, $\zeta$ is central (here $\zeta=Z^6$). Conjugacy classes
of reflections are powers of $c$. 

Group $G_9$ (order 192) is generated by $a=a_*,b=b_*Z^{-1},c=c_*$. 
The relations are $a^2=1$, $b^3=\zeta^{-1}$, $c^4=1$,
$abc=1$, $\zeta$ is central (here $\zeta=Z^3$). Conjugacy classes
of reflections are powers of $a,c$. 

Group $G_{10}$ (order 288) is generated by $a=a_*Z^{-1},b=b_*,c=c_*$. 
The relations are $a^2=\zeta^{-1}$, $b^3=1$, $c^4=1$,
$abc=1$, $\zeta$ is central (here $\zeta=Z^2$). Conjugacy classes
of reflections are powers of $b,c$. 

Group $G_{12}$ (order 48) 
is generated by $a=a_*,b=b_*Z^{-4},c=c_*Z^3$. 
The relations are $a^2=1$, $b^3=\zeta^{-1}$, $c^4=\zeta$,
$abc=1$, $\zeta$ is central (here $\zeta=Z^{12}$). Conjugacy classes
of reflections are powers of $a$. 

Group $G_{13}$ (order 96) is generated by
$a=a_*,b=b_*Z^2,c=c_*Z^{-3}$, and $f=c_*^2$. 
The relations are $a^2=1$, $b^3=\zeta$, 
$f^2=1$, $c^2\zeta=f$,
$abc=1$, $\zeta$ is central (here $\zeta=Z^6$). Conjugacy classes
of reflections are powers of $a$ and $f$.  

Group $G_{14}$ (order 144) 
is generated by $a=a_*,b=b_*,c=c_*Z^{-1}$. 
The relations are $a^2=1$, $b^3=1$, $c^4=\zeta^{-1}$,
$abc=1$, $\zeta$ is central (here $\zeta=Z^4$). Conjugacy classes
of reflections are powers of $a,b$. 

Group $G_{15}$ (order 288) 
is generated by $a=a_*,b=b_*,c=c_*Z^{-1},f=c_*^2$. 
The relations are $a^2=1$, $b^3=1$, $f^2=1$, 
$c^2\zeta=f$, 
$abc=1$, $\zeta$ is central (here $\zeta=Z^2$). Conjugacy classes
of reflections are powers of $a,b$, and $f$. 

{\bf Icosahedral family.}

The icosahedral group $G$: $a^2=b^3=c^5=1$, $abc=1$ (order 60).

The maximal group $\widehat G=G_{19}$ (order 3600):
$a_*^2=b_*^3=c_*^5=1$, $a_*b_*c_*=Z$, $Z$ is central; conjugacy classes
of reflections are powers of $a_*,b_*,c_*$. 

Group $G_{16}$ (order 600) is generated by $a=a_*Z^{-3},b=b_*Z^2,c=c_*$. 
The relations are $a^2=\zeta^{-1}$, $b^3=\zeta$, $c^5=1$,
$abc=1$, $\zeta$ is central (here $\zeta=Z^6$). Conjugacy classes
of reflections are powers of $c$. 

Group $G_{17}$ (order 1200) is generated by $a=a_*,b=b_*Z^{-1},c=c_*$. 
The relations are $a^2=1$, $b^3=\zeta^{-1}$, $c^5=1$,
$abc=1$, $\zeta$ is central (here $\zeta=Z^3$). Conjugacy classes
of reflections are powers of $a,c$. 

Group $G_{18}$ (order 1800) is generated by $a=a_*Z^{-1},b=b_*,c=c_*$. 
The relations are $a^2=\zeta^{-1}$, $b^3=1$, $c^5=1$,
$abc=1$, $\zeta$ is central (here $\zeta=Z^2$). Conjugacy classes
of reflections are powers of $b,c$. 

Group $G_{20}$ (order 360) 
is generated by $a=a_*Z^{-5},b=b_*,c=c_*Z^4$. 
The relations are $a^2=\zeta^{-1}$, $b^3=1$, $c^5=\zeta^2$,
$abc=1$, $\zeta$ is central (here $\zeta=Z^{10}$). Conjugacy classes
of reflections are powers of $b$. 

Group $G_{21}$ (order 720) is generated by $a=a_*,b=b_*,c=c_*Z^{-1}$. 
The relations are $a^2=1$, $b^3=1$, $c^5=\zeta^{-1}$,
$abc=1$, $\zeta$ is central (here $\zeta=Z^5$). Conjugacy classes
of reflections are powers of $a,b$. 

Group $G_{22}$ (order 240) 
is generated by $a=a_*,b=b_*Z^{5},c=c_*Z^{-6}$. 
The relations are $a^2=1$, $b^3=\zeta$, $c^5=\zeta^{-2}$,
$abc=1$, $\zeta$ is central (here $\zeta=Z^{15}$). Conjugacy classes
of reflections are powers of $a$. 

Note that in each case, conjugacy classes of reflections 
are represented by those of the generators $g$ for which one of the
defining relations is $g^p=1$. The braid groups $B_\Gamma$ 
of $\Gamma$ are obtained by removing such relations, and thus the
Hecke algebras are obtained by replacing them with the relations \linebreak
$(g-\bold b_1)...(g-\bold b_p)=0$, where
$\bold b_1,...,\bold b_p$ are invertible parameters (i.e., characters
of $\bold T$); see \cite{BMR,
Ma}. 

Now we proceed to prove the theorem. 
Since it is known that generically $\dim H(\Lambda)\ge |\Gamma|$, 
it suffices to show that $H$ is a finitely generated module over
$\CC[\bold T]$.

To do so, we note that 
the group $G$ is the group of even elements in a finite Coxeter group 
of rank 3: of type $A_1\times I_m$ in the dihedral case 
(which we are not considering) and $A_3,B_3,H_3$ for the tetrahedral, octahedral,
and icosahedral cases, respectively. Therefore, by Theorem 2.3 of \cite{ER}, 
$H$ is a finite module over $\CC[\bold T][Z,Z^{-1}]$ 
generated by $|G|$ generators (denoted in \cite{ER} by $T_w$,
$w\in G$). Let $S$ be the support of this module 
(a subvariety in $\bold T\times \CC^*$). If $(\Lambda,z)\in S$, 
then $H(\Lambda)$ must have an irreducible representation of dimension 
$d\le |G|^{1/2}$, in which $Z$ acts by the scalar $z$. 
Taking the determinants of the defining relations of $H$ in this
representation, we obtain that $z^d$ is a certain character of
$\bold T$ evaluated at $\Lambda$ 
(there are finitely many possibilities for such
characters). Thus $S$ is contained in a
finite union of subtori of $\bold T\times \CC^*$
whose projections to $\bold T$ are finite. We conclude that  
the action of $\CC[\bold T][Z,Z^{-1}]$
in $H$ factors through the algebra 
${\mathcal R}:=\CC[\bold T][Z,Z^{-1}]/J$ for 
a certain ideal $J$ whose zero set is $S$, and
the algebra ${\mathcal R}$ is a finitely generated module
over $\CC[\bold T]$. This implies the theorem.   

\subsection{The connection between $B(\lambda)$ and $H(\Lambda)$}

Now assume that $\Gamma=\widehat G$ is a maximal group, 
i.e. $G_7,G_{11},G_{19}$. In this case, using the notation of
Subsection \ref{nv}, the cyclotomic 
Hecke algebra is generated by elements $Y_k$,  
and invertible central variables $\bold b_{jk}$, $j=1,...,d_k$,
$k=1,...,m$, with defining relations 
$$
\prod_{j=1}^{d_k}  (Y_k-\bold b_{jk})=0, k=1,...,m,
$$
and $\prod_{k=1}^m Y_k=Z$ is central. 
The specialization $H(\Lambda)$ of the algebra $H$
is obtained when the variables $\bold b_{jk}$ map  
to complex numbers $\Lambda_{jk}$.  
 
Let $\Lambda=e^{2\pi i\lambda}$. In this case we can define an algebra
homomorphism $\phi: H(\Lambda)\to B(\lambda)$ 
(the Riemann-Hilbert homomorphism) in the following manner. 
Let $\zeta_k,k=0,1,...,m$ be distinct real numbers listed in
increasing order. Consider the 
differential equation 
$$
\frac{dF}{d\zeta}=\sum_{k=1}^m\frac{U_kF}{\zeta-\zeta_k}
$$
with respect to a holomorphic function $F(\zeta)$ with values in
$B(\lambda)$, and define $\phi(Y_k)$ to be the monodromy operator
of this equation around a loop
$\gamma_k$ which starts and ends at $\zeta_0$ and goes around
$\zeta_k$ counterclockwise, passing $\zeta_1,...,\zeta_{k-1}$
from below. It is easy to check that this gives rise to a well defined
homomorphism $\phi$. 

For any $s\in \Bbb C$, let $H(\Lambda)_s$ denote the 
generalized eigenspace of $Z$ with eigenvalue $e^{2\pi is}$. 
It is a direct summand subalgebra of $H(\Lambda)$. 
Similarly, denote by $B(\lambda)_s$ the 
direct summand subalgebra in $B(\lambda)$
which is the generalized eigenspace of $z$ with
eigenvalue $s$. It is clear that for any $s$, 
$\phi$ maps $H(\Lambda)_s$ to $B(\lambda)_s$. 

Also, let $H_*:=H(1)/(Z-1)$, and $B_0=B(0)/(z)=e_p\Pi_0e_p$. 
Since $\phi(Z)=e^{2\pi iz}$, the map $\phi$ descends to a
homomorphism $\phi_*:
H_*\to B_0$. 

\begin{proposition}\label{phi0}
(i) $\phi_*$ is an isomorphism. 

(ii) $\phi_0: H(1)_0\to B(0)$ 
is an isomorphism. In particular, 
the block of the trivial representation of $H(1)$ 
(i.e. the 1-dimensional representation where all the generators
act by 1) is equivalent to the category of representations 
of $B(0)$. 

(iii) There exists $\varepsilon>0$ such that if
$|\lambda_{jk}|<\varepsilon$ for all $j,k$ and $|s|<\varepsilon$, then
$\phi$ defines an isomorphism $H(\Lambda)_s\to B(\lambda)_s$. 
\end{proposition}

\begin{proof}
It is clear that $\phi(Y_k-1)=2\pi iU_k+h.d.t.$, 
where $h.d.t.$ denotes the higher degree terms in $U_k$.  
Thus the homomorphism $\phi: H(1)\to B(0)$ is surjective,
and hence $\phi_*$ is surjective. 
On the other hand, it is clear that $\phi(Z)=e^{2\pi iz}$, 
so $\phi$ factors through $H(1)_0$ where it descends to
$\phi_0$. So $\phi_0$ is surjective.

Now let us show that $\phi_*$ is an isomorphism. 
To do so, recall that the dimension of $B_0$ is $|G|$ 
(see e.g. \cite{MOV}). Thus it suffices to show that 
$\dim H_*\le |G|$. But this is trivial for type $A$ and 
follows from \cite{ER}, 
Theorem 2.3, applied to finite Coxeter groups of rank 3, for
types D and E. So (i) is proved. 

To prove (ii), it remains to show that $\phi_0$ is an isomorphism. 
For this purpose note that by Schur's lemma, any irreducible representation 
of $H(1)_0$ factors through $H_*$ and hence by (i) is 1-dimensional, with
all generators acting by $1$. This means that 
$H(1)_0={\rm limproj}_{n\to \infty}H(1)/J^n$, where 
$J$ is the ideal generated by the elements $Y_k-1$. 
So we can define elements $U_k'=(2\pi i)^{-1}\log(Y_k)\in H(1)_0$
(using the power series for the logarithm). The elements 
$U_k'$ satisfy the equation $(U_k')^{d_k}=0$, and 
$\sum_k U_k'+h.d.t.=0$, where $h.d.t.$ stand for higher degree
expression in $U_k'$. Since $U_k'$ obviously generate $H(1)_0$, 
we can put on $H(1)_0$ a decreasing filtration defined by
$\text{deg}(U_k')=1$. Then we find that ${\rm gr}H(1)_0$ is a quotient of
$B(0)$, which implies that $\dim H(1)_0\le \dim B(0)$. 
This implies (ii).

Finally, for proof of (iii) it suffices to note that 

1) $B(\lambda)$ is a flat deformation of $B(0)$, so $\dim
B(\lambda)=\dim B(0)$; 

2) There exists a constant $K>0$ such that for any 
$\varepsilon>0$, if $|\lambda_{jk}|<\varepsilon$ for all $j,k$, 
then $B(\lambda)=\oplus_{s: |s|<K\varepsilon}B(\lambda)_s$, 

3) The map $\phi_s: H(\Lambda)_s\to B(\lambda)_s$ is surjective; and

4) by Theorem \ref{finitel}, there is $\varepsilon>0$ such that 
if $|\lambda_{jk}|<\varepsilon$ for all $j,k$, then the 
dimension of the direct sum of $H(\Lambda)_s$ over $s$ such that 
$|s|<K\varepsilon$ is $\le \dim H(1)_0=\dim B(0)$. 
\end{proof}


\begin{thebibliography}{999}

\bibitem[BaKi]{BaKi}
Bakalov, B., Kirillov, A., Jr.,
 Lectures on tensor categories and modular functors. 
University Lecture Series, 21. American Mathematical Society, Providence, RI, 2001.

\bibitem[BMR]{BMR}
M. Brou\'e, G. Malle and R. Rouquier, Complex reflection groups,
braid groups, Hecke algebras, J. Reine Angew. Math. 500 (1998),
127-190.

\bibitem[CB1]{CB1} W. Crawley-Boevey, {\em
 Geometry of the moment map for representations of quivers.}
  Compositio Math. {\bf   126} (2001),   257--293.

\bibitem[CB2]{CB2} W. Crawley-Boevey, Annihilating elements
in generic deformed preprojective algebras, preprint. 

\bibitem[CBH]{CBH} W. Crawley-Boevey, M. Holland: {\it
 Noncommutative deformations of Kleinian singularities.} Duke
Math. J. {\bf 92}
(1998), 605--635.

\bibitem[EG]{EG} 
Etingof, P.; Ginzburg, V., 
Symplectic reflection algebras, Calogero-Moser space, and
deformed Harish-Chandra homomorphism.  Invent. Math.  147  (2002),  no. 2, 243--348.

\bibitem[EO]{EO}
Etingof, P., Ostrik, V., 
Module categories over representations of ${\rm SL}\sb q(2)$ and
graphs. Math. Res. Lett.  11  (2004),  no. 1, 103--114.

\bibitem[ER]{ER} P. Etingof and E. Rains, 
New deformations of group algebras of Coxeter groups, math. QA/0409261.


\bibitem[GGOR]{GGOR} Ginzburg, V.; Guay, N.; Opdam, E.; Rouquier, R.
On the category $\mathcal O$ for rational Cherednik algebras.
Invent. Math.  154  (2003),  no. 3, 617--651.

\bibitem[GP]{GP} I.M. Gelfand and V.A. Ponomarev, Model algebras and
representations of graphs, Func. Anal. and Applic. v.13 (1979), no. 3, 1--12.

\bibitem[Ko]{Ko} B. Kostant, The Coxeter element and the
branching law for the finite subgroups of SU(2),
math.RT/0411142.

\bibitem[MOV]{MOV} A. Malkin, V. Ostrik, M. Vybornov,
{\em  Quiver varieties and Lusztig's algebra.}\hfill\break
{\tt{arXiv:math.RT/0403222}}.

\bibitem[Ma]{Ma} G. Malle,
Degr\'es relatifs des algebres cyclotomiques associ\'ees 
aux groupes de r\'eflexions
complexes de dimension deux, Finite Reductive Groups: Related
Structures and Representations (M. Cabanes, ed.), Progress in
Mathematics, vol. 141, Birkh\"auser, 1997, pp. 
311-332. 

\bibitem[Me]{Me} Anton Mellit,
Algebras Generated by Elements with Given Spectrum and 
Scalar Sum and Kleinian Singularities,
math.RA/0406119. 

\bibitem[Ru]{Ru} 
W. Rump, Doubling a path algebra, or how to extend indecomposable
modules to simple modules, Representation theory of groups,
algebras and orders (Costanta, 1995), An. Stiint. Univ. 
Ovidius Constanta Ser Mat. 4 2 (1996), 174-185.

\bibitem[ST]{ST} 
G. C. Shephard and J. A. Todd, Finite unitary reflection groups,
Canad. J. Math. v.6 (1954), p. 274=-304.

\bibitem[YZ]{YZ}
Yekutieli, A., Zhang, J. J.,
Rings with Auslander dualizing complexes.  J. Algebra  213  (1999),  no. 1, 1--51.

\end{thebibliography}
\end{document}